\documentclass[11pt]{article}
\usepackage{amssymb,amsmath,latexsym,amsthm} %
\usepackage{xspace,bm}
\usepackage{fullpage}
\usepackage{hyperref} %
\usepackage{natbib}  %
\usepackage{comment}

\hypersetup{colorlinks=true, citecolor=blue,
  bookmarksopen=false, pdfstartview=FitH, pdfview=FitH}
\bibpunct{[}{]}{,}{n}{,}{,}
\newtheorem{theorem}{Theorem}
\newtheorem{proposition}[theorem]{Proposition}
\newtheorem{corollary}[theorem]{Corollary}
\newtheorem{lemma}[theorem]{Lemma}

\theoremstyle{definition}
\newtheorem{definition}[theorem]{Definition}
\newtheorem{question}{Question}
\newtheorem{example}{Example}
\newtheorem{remark}[theorem]{Remark}

\newcommand{\Tau}{\mathcal{T}}

\DeclareMathOperator{\Var}{\mathrm{Var}}
\DeclareMathOperator{\Cov}{\mathrm{Cov}}
\DeclareMathOperator{\Ent}{\mathrm{Ent}}

\newcommand{\xxx}{\mathcal{X}} \newcommand{\yyy}{\mathcal{Y}}

\newcommand{\EEE}{\mathcal{E}} 

\newcommand{\nnn}{\mathbb{N}}
\newcommand{\ppp}{\mathbb{P}} 
\newcommand{\qqq}{\mathbb{Q}} \newcommand{\rrr}{\mathbb{R}}

 \newcommand{\EXP}{\mathbb{E}}
\newcommand{\zzz}{\mathbb{Z}} 
\newcommand{\sep}{\rm{\bf sep}}
\newcommand{\D}{{\rm {\bf D}}}

\newcommand\secref[1]{\hyperref[#1]{Section \ref*{#1}}} %
\newcommand\lemref[1]{\hyperref[#1]{Lemma \ref*{#1}}}
\newcommand\corref[1]{\hyperref[#1]{Corollary \ref*{#1}}}
\newcommand\egref[1]{\hyperref[#1]{Example \ref*{#1}}}
\let\eqnref\autoref

\newcommand{\rv}[1]{\mathbf{#1}}

\title{A Tight Bound for the Lamplighter Problem}
\author{Murali K. Ganapathy\thanks{Google, 1600 Amphitheatre Pkwy, Mountain
View, CA - 94043} \and Prasad Tetali\thanks{School of Mathematics
and College of Computing, Georgia Tech., Atlanta, GA 30332-0160.
Research supported in part by NSF grant DMS-0401239.}}
\date{}
\begin{document}

\maketitle

\begin{abstract}
We settle an open problem, raised by Y. Peres and D. Revelle,
concerning the $L^2$ mixing time of the random walk on the
lamplighter graph. We also provide general bounds relating the
entropy decay of a Markov chain to the separation distance of the
chain, and show that the lamplighter graphs once again provide
examples of tightness of our results.
\end{abstract}

\section{Introduction}

Given a finite connected graph $G$, a vertex of the lamplighter
graph $G^\lozenge$ consists of a 0-1 labeling of the vertices of
$G$, and a marked vertex of $G$. Each vertex has a lamp, the marked
vertex indicates the position of a lamplighter and the labeling at
any time indicates the off-on status of each lamp (vertex). The
lamplighter random walk on $G^\lozenge$ corresponds to the
lamplighter performing a random walk on $G$, while randomizing the
status of each lamp, as he/she visits the corresponding vertex. When
$G$ is a cycle or a complete graph, the corresponding lamplighter
chains were studied by H\"aggstr\"om and Jonasson \cite{HJ}.
Vertex transitive, other special classes and more general  graphs
were considered in detail by Peres and Revelle
\cite{LampLighter}, who provided general upper and lower bounds for
mixing times of the lamplighter random walk.

By tightening the analysis in \cite{LampLighter}, we prove an
optimal upper bound on the $L^2$ mixing time of the lamplighter
Markov chain on  a class of graphs considered in \cite{LampLighter}. The
mixing time and related measures are defined via

\begin{definition}
Let $G$ be a connected $d$-regular undirected graph.
\begin{itemize}
\item The relaxation time
$T_{rel}(G)=\max_\lambda 1/(1-|\lambda|)$ where the maximum is taken over
non-trivial eigenvalues $\lambda$ of the normalized adjacency matrix of $G$.
\item $\tau(G)$ ($\tau_2(G)$) is the time for the random walk on $G$ to be
within $1/4$ of uniform distribution in total variation distance
($L_2$-distance, respectively).
\item $H(G)$ is expected time it takes for the random walk to travel
from $x$ to $y$, where the choice of $x$ and $y$ is adversarial.
\end{itemize}
\end{definition}

There is a more popular definition of relaxation time  (see e.g.
\cite{AF, MT}). We choose this definition as it is easier to work
with. Standard inequalities imply that the two numbers differ by
at most 1.

For a graph $G$, let
$\tau_2(G^{\lozenge})$ denote the $L^2$ mixing time of the
lamplighter random walk on $G^{\lozenge}$.  Then  our main theorem
is as follows.

\begin{theorem}
\label{main_thm} Suppose $G$ is a regular undirected graph for which, $H(G) \le
\kappa |G|$, for some universal constant $\kappa>0$. Then there
exists a constant $c=c(\kappa)$ such that
\begin{equation}
\label{main_ub}
\tau_2(G^{\lozenge}) \le c \, |G|\bigl(T_{rel}(G) + \log|G|  \bigr)\, .
\end{equation}
\end{theorem}

The above theorem refines and improves upon (by providing a matching
upper bound to) the result of Y. Peres and D. Revelle
\cite{LampLighter},  who proved  (inter alia) the following theorem.

\begin{theorem}[Peres-Revelle]
Under the hypothesis of Theorem~\ref{main_thm}, there exist
constants, $c_1, c_2$ depending on $\kappa$ such that
\begin{equation}
\label{pr_bound}
c_1 \,  |G|\bigl(T_{rel}(G) + \log|G|  \bigr)\, \le
\tau_2(G^{\lozenge}) \le c_2 \, |G|\bigl(T_{\rm tv}(G) + \log|G|  \bigr)\,
,
\end{equation}
where $T_{\rm tv}(G)$ is the total variation mixing time of the
simple random walk on $G$.
\end{theorem}

Note that our theorem shows that for $G=\zzz_2^n$,  the correct
order of magnitude for $\tau_2(G^\lozenge)$ to be $n2^n$, since the
relaxation time is of order $n$. This settles
 Problem 4 mentioned at the end of \cite{LampLighter}, while our
theorem itself settles the question raised as Problem 5 in the
affirmative.

In \cite{LampLighter}, the lamplighter random walk on the two-dimensional
torus was shown to be an example of a chain for which the relaxation time,
the total variation mixing time, and the $L^2$ mixing time were all shown
to be distinct orders of magnitude.  In this paper, we use the one and two-
dimensional tori as examples which further separate the mixing time in entropy
(relative to stationarity) from the rest of the above mixing times.  These
examples also illustrate tightness of the following other result of this
paper. We show that in general the {\em entropy mixing time} is at worst a
factor of $\log \log(1/\pi_*)$
larger than the total variation mixing time for reversible Markov chains
(see Corollary~\ref{cor:ent_tv}) below.)
 This is accomplished by relating the relative entropy
 to  the so-called separation distance of a Markov chain.

\section{The Lamplighter Result}
In this section we derive some preliminary technical lemmas and a
key theorem from which the main theorem follows. Since random walks on regular
undirected graphs are equivalent to reversible Markov chain with uniform
stationary distribution, we will use the latter from now on.

Assume that $\ppp$ is a Markov chain with uniform stationary
distribution $\pi$ on a finite state space $\xxx$. Let $H =
\max_{x,y} \EXP_xT_y$ denote the  maximal hitting time (also
called the maximum expected first passage time) of the chain. Let
$T_{rel}$
%%= {1}/{\lambda}$
denote the relaxation time of the chain, where
$\lambda$ is the spectral gap of the chain.

\medskip

As observed by Peres and Revelle, the $L^2$ mixing time of the
lamplighter graph $G^{\lozenge}$ depends upon the moment generating
function of the cover time of the underlying graph $G$. More
precisely, if $S_t$ denotes the set of unvisited vertices (by the
lamplighter) by time $t$, then to get convergence in the $L^2$ (or
equivalently, in the uniform metric), one needs $\EXP 2^{|S_t|} \le
1+\epsilon$, for $\epsilon >0$. Our main technical contribution is
as follows.

\medskip

 Let $\ppp$ be a reversible Markov Chain on the state space $\xxx$
with $\pi$ as the stationary distribution.
\begin{theorem}\label{key_thm}
Let the chain given by $\ppp$ start in an initial
distribution $\mu$ so that $\mu \geq \pi/2$.  Let the maximal
hitting time $H=H(\ppp)$  satisfy $H \leq c_1|\xxx|$ for a constant
$c_1 \geq 1$.
%%$|\xxx| \geq 2$.
Let $\theta \geq 2$ be arbitrary, and let $\rv{S}_t$ denote the
set of vertices which have not been visited by time $t$.  Then
there exists a universal constant $c$ such that for all $a,b > 0$,
and  for $t \geq t'= C_1\, |\xxx| T_{rel}\log\theta + C_2\,
|\xxx|\log|\xxx|$, we have
% $t \geq t'=c\, c_1^2\, |\xxx|\bigl[2(1+a)T_{rel}\log\theta +
%(1+b)\log|\xxx|\big]$ we have
\begin{equation*}
\EXP[\theta^{|\rv{S}_t|}] \leq 1 + \delta + \delta^2 + \delta^9\,,
\end{equation*}
where  $\delta=\theta^{-(1+2a)T_{rel}}|\xxx|^{-b}$,
$C_1=2c\,c_1^2(1+a)$ and $C_2=c\,c_1^2(1+b)$. \noindent In
particular, when $a=b=1$, we have $\EXP[\theta^{|\rv{S}_t|}] <
1.21$.
 \end{theorem}

Once we have Theorem~\ref{key_thm}, the main theorem
(Theorem~\ref{main_thm}) follows in a straightforward way, as in
\cite{LampLighter}.
 We begin with a few simple lemmas.

Let $\ppp$ be a Markov chain with uniform stationary distribution
$\pi$.  Let $\sigma_1$ denote the second largest singular value of
$\ppp$. In particular, if $\ppp$ is reversible, then $\sigma_1 =
\max(\lambda_1,|\lambda_{N-1}|)$, where $1= \lambda_0 \ge \lambda_1
\ge \cdots \ge \lambda_{N-1}$ denote all the eigenvalues of $\ppp$.
Recall the following basic fact, whose proof we include for
completeness.

 \begin{lemma}\label{lem:fact}
 Let $f: \xxx \to \rrr$.
  Let $\rv{X}_i$ denote the state of the chain $\ppp$ at time $i$. Then
  \begin{equation*}
  \Cov(f(\rv{X}_1),f(\rv{X}_t)) \leq \sigma_1^t \Var(f(X_1))\, ,
  \end{equation*}
 where $\sigma_1$ is the second largest singular value of $\ppp$.
 \end{lemma}
\begin{proof}
Let $g=\ppp^tf$, so that $\Cov(f(\rv{X}_1),f(\rv{X}_t))=\Cov(f,g)$.
Hence we have
\begin{align*}
\Cov(f,g)
   &= \EXP_\pi [(f-Ef)(g-Eg)] \\
  &= \EXP_\pi [(f-Ef)(\ppp^tf-E\ppp^tf)] \\
  &= \EXP_\pi [(f-Ef)((\ppp^t-E)f)] && (E\ppp=E)\\
  &= \EXP_\pi [(f-Ef)((\ppp-E)^tf)] && \text{$\ppp^t-E=(\ppp-E)^t$}\\
  &\leq \sqrt{\EXP_\pi((I-E)f)^2} \, \sigma_1^t \, \sqrt{\Var_\pi(f)} \\
  &= \lambda_*^t\Var_\pi(f)
\end{align*}
where we used the fact that the operator norm of $\ppp-E$ is
$\sigma_1$.
\end{proof}

\begin{lemma}\label{lem:key}
  Let $\{X_t\}$ be a reversible Markov chain on $\xxx$, with uniform
  stationary distribution $\pi$. Assume that $\Pr\{X_0 = x \} \geq
  \pi(x)/2$. Let $T_{rel}$ be the relaxation time of the chain. Let
  $T^{+}_x$ denote the return time to $x$ and assume there are
  $\epsilon,\delta > 0$ for which
  \begin{equation*}
    \Pr_x (T^{+}_x \geq \epsilon |\xxx|) \geq \delta > 0
  \end{equation*}
  for all $x \in \xxx$. Let $\yyy \subset \xxx$ be such that $|\yyy|
  \geq T_{rel}$. Then the probability of hitting at least
  $\delta\epsilon |\yyy|/4$ elements of $\yyy$ by time
  $ C\delta^{-2}T_{rel}/\pi(\yyy)$ is at least $1/2$,
  where $C \geq 16$ is an absolute constant.
\end{lemma}
\begin{proof}
  Let $r=\epsilon\, T_{rel}/\pi(\yyy)$. For $1 \leq i
  \leq r$, let $I_i$ be an indicator random variable for the event
  $\{X_i \in \yyy\}$ and $J_i$ for the event $\{X_i \in \yyy\}$ and
  $\{X_j \neq X_i\}$ for $i < j \leq r$. Finally let $J = \sum_i
  J_i$ and $I=\sum_i I_i$.

\

  Note that $J$ is the number of distinct elements of $\yyy$ which
  have been visited in the time interval $[1,r]$. Also we have $\Pr
  \{ J_i = 1 | I_i = 1 \} \geq \delta$ since $r \leq
  \epsilon|\xxx|$. This together with the fact that $\EXP [I_i] \geq
  \pi(\yyy)/2$ (due to our assumption on the initial distribution), gives
  \begin{equation}\label{eqn:A}
    \EXP [J] \geq \delta \EXP [I] \geq \delta r \pi(\yyy)/2 =
    \frac{\delta\epsilon T_{rel}}{2}\, .
  \end{equation}

  To conclude $\Pr \{ J \geq \epsilon\delta T_{rel}/4 \}$ is bounded
  away from 0, we bound $\EXP[I^2]$.

  Fix $1 \leq i \leq j \leq r$. From \lemref{lem:fact} we have
  \begin{equation*}
    \sum_{j \geq i} \Cov(I_i,I_j) \leq \sum_{j \geq i}
\lambda_*^{j-i} \Var(I_i) \leq \frac{\Var(I_i)}{1-\lambda_*} =
T_{rel}\Var(I_i) \leq T_{rel} \EXP[I_i]\, .
\end{equation*}
  Now we have
  \begin{equation*}
  \EXP[I^2] \leq 2 \sum_i \sum_{j \geq i} \Cov(I_i,I_j) \leq 2
  \sum_i T_{rel} \EXP[I_i] = 2 T_{rel} \EXP[I]\, .
  \end{equation*}

  Since $J_i \leq I_i$ we have $\EXP[J^2] \leq \EXP[I^2] \leq 2
  T_{rel} \EXP[I] \leq 2 T_{rel} \EXP[J] / \delta$, since $\EXP[J_i]
  \geq \delta \EXP[I_i]$. Hence using \eqnref{eqn:A} we have
  \begin{equation}\label{eqn:B}
  \EXP[J^2] \leq \frac{4}{\delta^2\epsilon}\EXP[J]^2\, .
  \end{equation}

  Now let $\alpha$ be the indicator for the event $J \geq \EXP[J]/2$.
  Then $\EXP[J(1-\alpha)] \leq \EXP[J]/2$ and hence $\EXP[J\alpha] \geq
  \EXP[J]/2$. Now by Cauchy-Schwartz, we have $\EXP[J\alpha]^2 \leq
\EXP[J^2]\EXP[\alpha^2]$. Hence
  \begin{equation*}
  \Pr \{ J \geq \EXP[J]/2 \} = \EXP[\alpha^2]
  \geq \frac{\EXP[J\alpha]^2}{\EXP[J^2]}
  \geq \frac{(\EXP[J]/2)^2}{\EXP[J^2]}
  \geq \frac{\delta^2\epsilon}{16}\, .
  \end{equation*}
  Thus in a trial of length $r=\epsilon |\xxx|T_{rel}/|\yyy|$, the
  probability that we do not pick up $\delta\epsilon|\yyy|/4$
  elements of $\yyy$ is less than $1 - \delta^2\epsilon/16$. Hence
  if we repeat this for $C\delta^{-2}/\epsilon$ intervals of length
  $r$ each, we can reduce the probability of failure to less than
  $1/2$. Note that $C \geq 16$ here.
\end{proof}

Suppose that the initial distribution $\mu$ is such that $\mu  \geq
\pi/2$. Then if the current set of unvisited states is large ($\geq
T_{rel}$) then \lemref{lem:key} shows that we visit
$\Omega(T_{rel})$ new states within time $O(T_{rel}/\pi(\yyy))$ with
probability $\geq 1/2$. Once the set of unvisited states gets
smaller than $T_{rel}$ things are in better shape. The next lemma
establishes the assumption of \lemref{lem:key} and handles the case
when the set of unvisited vertices is small.

\begin{lemma}\label{lem:excursion}
  Let $\ppp$ be a Markov chain with uniform stationary distribution
  and maximal hitting time $H$. For $x \in \xxx$, let $T^{+}_x$
  denote the expected length of the return time to $x$. Then
  \begin{equation*}
  \min_x \Pr_x \Bigl\{ T^{+}_x \geq \frac{\xxx|}{2} \Bigr\} \geq
  \frac{|\xxx|}{2H}\, .
  \end{equation*}
  Also, for any $\yyy \subseteq \xxx$, with probability $\geq 1/2$,
  we visit at least $|\yyy|/2$ elements of $\yyy$ by time $4H$.
\end{lemma}
\begin{proof}
  Since the stationary distribution is uniform, $\EXP_x[T^{+}_x] =
  |\xxx|$. If after $|\xxx|/2$ steps we have not yet returned
  to $x$, and are currently at state $y$, then we expect to visit
  $x$ within another $H$ steps. Hence
  \begin{equation*}
    |\xxx| = \EXP_x[T^{+}_x] \leq \Pr\{ T^{+}_x \leq |\xxx|/2 \}
    |\xxx|/2 + \Pr\{ T^{\gamma} \geq |\xxx|/2 \} H\, .
  \end{equation*}
  Rearranging terms, we get the result.

  For the second result:
  Fix $x \in \yyy$ and let $\rv{H}_x$ denote the time when
  $x$ is visited. $\EXP[\rv{H}_x] \leq H$. Thus by time $4H$ we visit
  $x$ with probability $\geq 3/4$.

  Let $Y$ denote the number of elements of $\yyy$ which have been
  visited by time $4H$. Then $\EXP[Y] \geq 3|\yyy|/4$. If $q=\Pr\{Y
  \geq |\yyy|/2\}$, we have
  \begin{equation*}
    3|\yyy|/4 \leq \EXP[Y] \leq (1-q)\frac{|\yyy|}{2} + q|\yyy|\, .
  \end{equation*}
  Solving for $q$  gives $q \geq 1/2$.
\end{proof}

%\begin{theorem}\label{thm:mainresult}
%Let $\ppp$ be a reversible Markov Chain with initial distribution
%$\mu \geq \pi/2$. Let $H \leq K|\xxx|$ for a constant $K \geq 1$ and
%$|\xxx| \geq 2$. Let $\theta \geq 2$ be arbitrary. There exists a
%universal constant $c$ such that for all $a,b > 0$ and $t \geq
%t'=cK^2|\xxx|\bigl(2(1+a)T_{rel}\log\theta + (1+b)\log|\xxx|\bigr)$
%we have
%\begin{equation*}
%\EXP[\theta^{|\rv{S}_t|}] \leq 1 + \delta + \delta^2 + \delta^9
%\end{equation*}
%where $\rv{S}_t$ denotes the set of vertices which have not been
%visited at time $t$ and $\delta=\theta^{-(1+2a)T_{rel}}|\xxx|^{-b}$.
%In particular when $a=b=1$, we have $\EXP[\theta^{|\rv{S}_t|}] <
%1.21$.
%\end{theorem}
%\begin{proof}

\

\noindent
{\bf Proof of Theorem~\ref{key_thm}}.
Let $r=\lfloor|\xxx|/T_{rel}\rfloor$ and for $i=0,1,\dots,r-1$, let
$k_i = |\xxx|- iT_{rel}$ and for $i=r$, $k_r = 0$. Define stopping
times $\rv{T}_i$ as the time when $|\rv{S}_t| = k_i$ for the first
time.

\

\noindent
From \lemref{lem:excursion} and \lemref{lem:key}, it then follows
that $\rv{T}_i - \rv{T}_{i-1}$ is stochastically dominated by
$\rv{\alpha}_i = \gamma/k_i \rv{Z}_i$ where $\rv{Z}_i$ is geometric
with mean $2$ and $\gamma = C\cdot K^2|\xxx|T_{rel}$ and $C \geq 16$
is a universal constant.

\

\noindent
Fix $t > 0$, $i < r$ and $\beta > 0$ be arbitrary, Then
\begin{align*}
 \Pr \{ \rv{T}_i \geq t \}
&  = \Pr \Bigl\{ \sum_{j=1}^i \alpha_i \geq t \Bigr\}
 \leq \Pr \Bigl\{ \sum_{j=1}^i \frac{\gamma}{k_i}\rv{Z}_i
    \geq t \Bigr\} \\
 &\leq \exp(-t\beta) \EXP \Bigl[ \sum_{j=1}^i
 \frac{\gamma\beta}{k_i}\rv{Z}_i \Bigr]
 \leq \exp(-t\beta) \prod_{j=1}^i \EXP \Bigl[
 \frac{\gamma\beta}{k_i}\rv{Z}_i \Bigr]\, .
\end{align*}
Choose $\beta$ so that $\beta=k_i/3\gamma$ so that $\gamma\beta \leq
k_j/3$ for all $j \leq i$. For $\alpha \leq 1/3$, $\EXP[\alpha
\rv{Z}_j] \leq \exp(3\alpha)$. This gives
\begin{equation*}
\EXP \Bigl[ \frac{\gamma\beta}{k_j}\rv{Z}_j \Bigr] \leq
\exp(k_i/k_j)\, .
\end{equation*}
Hence we have
\begin{equation*}
 \Pr \{ \rv{T}_i \geq t \} \leq \exp\Bigl( -t\frac{k_i}{3\gamma} +
 \sum_{j=1}^i \frac{k_i}{k_j} \Bigr) \leq \exp\Bigl(
 -t\frac{k_i}{3\gamma} + \frac{k_i}{T_{rel}} \log|\xxx| \Bigr)\, .
\end{equation*}

\noindent For $i=r$, \lemref{lem:excursion} implies $(\rv{T}_r -
\rv{T}_{r-1})$ is stochastically dominated by the sum of
$\ell=\log_2(2T_{rel})$ independent geometric random variables
with mean $4K|\xxx|$. Applying a Chernoff bound we get
\begin{equation}
\Pr \{ \rv{T}_r - \rv{T}_{r-1} \geq t \}
 \leq \exp\Bigl(\ell-\frac{t}{4K|\xxx|}+
 \ell\log\Bigl(\frac{t}{4K|\xxx|}\Bigr)\Bigr)\, .
\end{equation}

\noindent
Breaking the values of $|S_t|$ into intervals of size $T_{rel}$ we
have
\begin{equation*}
\EXP[\theta^{|S_t|}]
 \leq 1 + \sum_{i=0}^r \theta^{k_i+T_{rel}}
      \Pr\{ |S_t| \geq k_i \}
 = 1 + \sum_{i=0}^r \theta^{k_i+T_{rel}}
 \Pr \{\rv{T}_i \geq t \}\, .
\end{equation*}

\noindent
For $i < r$, $k_i \geq T_{rel}$ and hence
\begin{equation}\label{eqn:butlast}
\begin{split}
  \theta^{k_i+T_{rel}}\Pr \{ \rv{T}_i \geq t \} &\leq
  \exp\Bigl( (k_i+T_{rel})\log \theta -t\frac{k_i}{3\gamma} +
       \frac{k_i}{T_{rel}} \log|\xxx| \Bigr) \\
  &\leq \exp\Bigl(2k_i\log\theta + \frac{k_i}{T_{rel}}\log|\xxx| - t
\frac{k_i}{3\gamma}\Bigr)\, .
\end{split}
\end{equation}

\noindent
When $i=r$, $0=k_r \leq T_{rel} \leq k_{r-1}$ and hence
\begin{equation}\label{eqn:last}
\begin{split}
\theta^{k_r+T_{rel}}\Pr \{ \rv{T}_r \geq t \}
  &\leq \theta^{T_{rel}}\bigl(
     \Pr \{ \rv{T}_{r-1} \geq t/2 \} +
     \Pr \{ \rv{T}_r - \rv{T}_{r-1} \geq t/2 \} \bigr) \\
  &\leq \theta^{T_{rel}} \exp\Bigl(
     -\frac{t}{2}\frac{k_{r-1}}{3\gamma} + \frac{k_{r-1}}{T_{rel}} \log|\xxx|
     \Bigr) \\
  &\quad + \theta^{T_{rel}}
     \exp\Bigl(\ell-\frac{t}{8K|\xxx|}+
 \ell\log\Bigl(\frac{t}{8K|\xxx|}\Bigr)\Bigr)\, ,
\end{split}
\end{equation}
where $\ell=\log_2(2T_{rel})$.

\

\noindent
Let $t' = 6CK^2|\xxx|\bigl(2(1+a)T_{rel}\log\theta+ (1+b)\log|\xxx|
\bigr)$ for any $a,b > 0$ and hence take $c=6C$. We now show that
for $t \geq t'$, $\EXP[\theta^{|\rv{S}_{t'}|}] -1$ is small.
Recall that $\gamma=CK^2|\xxx|T_{rel}$, hence we have
\begin{align*}
\frac{t'}{3\gamma} &= 4(1+a)\log \theta + 2(1+b)\frac{\log
|\xxx|}{T_{rel}} \\
\frac{t'}{8K|\xxx|} &= \frac{3CK}{4}\bigl(2(1+a)T_{rel}\log\theta +
(1+b)\log|\xxx|\bigr) \\
\ell' := \log\Bigl(\frac{t'}{8K|\xxx|}\Bigr) &\geq
\log\Bigl(\frac{3(1+a)CK\log\theta\cdot 2T_{rel}}{4}\Bigr) \geq
\log(6\cdot(2T_{rel}))\, ,
\end{align*}
since $C\geq 16, K\geq 1, \theta \geq 2$.

\

\noindent
Now for $t \geq t'$ and $i < r$, \eqnref{eqn:butlast} reduces to
\begin{equation}\label{eqn:final_butlast}
\begin{split}
\theta^{k_i+T_{rel}}\Pr \{ \rv{T}_i \geq t' \}
 &\leq \exp \Bigl(
   2k_i\log \theta + \frac{k_i}{T_{rel}}\log |\xxx|
   -4(1+a)k_i\log \theta -2(1+b)\frac{k_i}{T_{rel}}\log|\xxx|
   \Bigr) \\
 &= \theta^{-(2+4a)k_i}|\xxx|^{-(1+2b)k_i/T_{rel}}\, .
\end{split}
\end{equation}

\noindent
And \eqnref{eqn:last} reduces to
\begin{equation}\label{eqn:simple_last}
\begin{split}
\theta^{T_{rel}}\Pr \{ \rv{T}_i \geq t' \}
  &\leq \theta^{T_{rel}} \exp \Bigl( -2(1+a)k_{r-1}\log\theta
  -(1+b)\frac{k_{r-1}}{T_{rel}}\log{\xxx} +
  \frac{k_{r-1}}{T_{rel}}\log{\xxx} \Bigr) \\
  &\quad +\theta^{T_{rel}} \exp\Bigl(\ell-\exp(\ell')+
 \ell\ell'\Bigr) \\
 &\leq \exp\Bigl(
       -(1+2a)k_{r-1}\log\theta
       -b\frac{k_{r-1}}{T_{rel}}\log|\xxx| \Bigr) \\
 &\quad +
       \exp\Bigl( T_{rel}\log\theta + (1+\ell')(\ell'-\log(6))
       - \exp(\ell') \Bigr)\, ,
\end{split}
\end{equation}
using $\ell=\log_2(2T_{rel})$ and $\ell' \geq \ell + \log(6)$. Here
we also use the fact that $f(x)= \exp(x)/4 - (1+x)(x-\log(6)) \geq 0$ for
all $x \geq 0$, and $k_{r-1} \geq T_{rel}$.

\

\noindent
We now have
\begin{equation}\label{eqn:final_last}
\begin{split}
\theta^{T_{rel}}\Pr \{ \rv{T}_i \geq t' \}
 &\leq \theta^{-(1+2a)T_{rel}}|\xxx|^{-b} \\
 &\quad + \exp\Bigl( T_{rel}\log\theta -
          \frac{9CK}{16}\Bigl(2(1+a)T_{rel}\log\theta + (1+b)\log|\xxx|
    \Bigr)\Bigr) \\
 &\leq \theta^{-(1+2a)T_{rel}}|\xxx|^{-b}
   + \exp\Bigl(-9(1+2a)T_{rel}\log\theta - 9(1+b)\log|\xxx|\Bigr)
   \\
 &\leq \theta^{-(1+2a)T_{rel}}|\xxx|^{-b} +
   \theta^{-9(1+2a)T_{rel}}|\xxx|^{-9b}\, ,
\end{split}
\end{equation}
using $C \geq 16$.

\

\noindent
Combining \eqnref{eqn:final_last} and \eqnref{eqn:final_butlast} we
have
\begin{equation}\label{eqn:butfinal}
\EXP[\theta^{|\rv{S}_{t'}|}]
 \leq 1 + \sum_{i=0}^{r-1} \theta^{-(2+4a)k_i}|\xxx|^{-(1+2b)k_i/T_{rel}}
 + \theta^{-(1+2a)T_{rel}}|\xxx|^{-b} +
 \theta^{-9(1+2a)T_{rel}}|\xxx|^{-9b}\, \, .
\end{equation}

\

\noindent
Now let $\eta=\theta^{-(1+2a)T_{rel}}$. Using $k_{r-i} \geq
iT_{rel}$, we get
\begin{equation*}
\EXP[\theta^{|\rv{S}_{t'}|}] \leq 1+\sum_{i=0}^{r-1}
\eta^{2(r-i)}|\xxx|^{-(1+2b)(r-i)} + \eta|\xxx|^{-b} +
\eta^9|\xxx|^{-9b}\, \, .
\end{equation*}

\noindent
Summing up the geometric progression and simplifying we get
\begin{equation}\label{eqn:final}
\EXP[\theta^{|\rv{S}_{t'}|}]
  \leq 1 + \eta^2|\xxx|^{-2b} +
  \eta|\xxx|^{-b} + \eta^9|\xxx|^{-9b}\, \, .
\end{equation}

\noindent Observe that by setting $a=b=1$ and using $\theta\geq
2,T_{rel}\geq 1/2,|\xxx|\geq 2$, we get $\eta \leq 2^{-3/2}$ and
$|\xxx|^{-b} \leq 1/2$. Finally, substituting in \eqnref{eqn:final}
we have $\EXP[\theta^{|\rv{S}_{t'}|}] < 1.21\, . \hfill \Box $
%%\end{proof}

\

\noindent {\bf Proof of Theorem~\ref{main_thm}}. Implicit in the
proof of Theorem 1.4 of \cite{LampLighter} is that for a suitable constant
$C>0$,
\begin{equation}\label{lamp_eqn_14}
\tau_2(G^\lozenge) \leq C\left( \tau_2(G) + \min_t
\left\{\EXP\bigl[2^{|S_t|}\bigr] < 2\right\} \right)\,.
\end{equation}

Since $G$ is undirected, the random walk on $G$ is reversible.
Hence after time $4\tau(G)$, the distribution $\mu$ of the
lamplighter's position satisfies $\mu \geq \pi/2$. Since the
random walk on $G$ is reversible, we have (e.g., by \cite{Aldous})
that the mixing time is bounded above by maximal hitting time
$H(G)$. Thus by running the Lamplighter chain for an initial
$O(H(G))$ steps, we can ensure that the assumption of
Theorem~\ref{key_thm} holds.

Hence Theorem~\ref{key_thm} implies $\EXP\bigl[2^{|S_t|}\bigr] <
2$ for $t = O(H(G)) + O(|G|\cdot(T_{rel}+\log|G|))$. Thus
(\ref{lamp_eqn_14}) gives
\begin{equation}
\tau_2(G^\lozenge) = O\bigl( \tau_2(G) +
|G|\cdot(T_{rel}+\log|G|)\bigr),
\end{equation}
since $H(G)\leq\kappa|G|$. Finally, the regularity of $G$ implies
that the stationary distribution of the random walk on $G$ is
uniform. Thus $\tau_2(G) = O(\tau(G)\log|G|)$ implying that
$\tau_2(G) = O(|G|\log|G|). \hfill \Box$

\section{Separation Distance and Entropy Decay}
Recall that Pinsker's inequality lets one bound the total
variation mixing time of a Markov chain {\em from above} by the
entropy decay time, up to an absolute constant. In this section,
we show that for reversible Markov chains, the time for the
relative entropy to decay to within $1/e$ is no larger than $\log
\log(1/\pi_*)$ times that of the total variation mixing time. We
actually prove a more general result for all Markov chains from
which the above will follow under the additional assumption of
reversibility.

%%Furthermore we prove that our bound is essentially tight using a
%lamplighterchain as %%an example.

\

Once again let $\ppp$ be a Markov kernel with stationary
distribution $\pi$ on a (finite) state space $\xxx$. First recall
the definition of {\em separation} between a chain at time $n\ge 0$
and $\pi$.

\begin{definition}
For $n \in \nnn$ and $x \in \xxx$,  set
\begin{equation*}
\sep_{\ppp}(x,n) = \max_{y\in \xxx} \Bigl(1 -
\frac{\ppp^n(x,y)}{\pi(y)} \Bigr) \, , \ \ \sep_{\ppp}(n) = \max_{x
\in \xxx} \, \sep_{\ppp}(x,n) \, .
\end{equation*}

Also set
\begin{equation*}
d_{\rm tv}(\ppp^n(x, \cdot) - \pi) = \sum_{y \in \pi(y)>
\ppp^n(x,y)} (\pi(y) - \ppp^n(x,y))\,,
\end{equation*}
\begin{equation*}
\|\ppp^n(x, \cdot)/\pi(\cdot) - 1\|_2 = \Bigl(\sum_{y \in \xxx}
\bigl(\ppp^n(x,y)/\pi(y) -1\bigr)^2 \pi(y)\Bigr)^{1/2}\,.
\end{equation*}
\end{definition}

When understood from the context, we drop the subscript in
$\sep_{\ppp}$. Recall that the function $n \mapsto \sep(n)$ is
non-increasing and sub-multiplicative (see \cite{AD,AF} for more
details.)  It is well-known and is easily seen that $d_{\rm
tv}(\ppp^n(x, \cdot) - \pi) \le \sep(x,n)$:
\begin{eqnarray*}
d_{\rm tv}(\ppp^n(x, \cdot) - \pi)
& = & \sum_{y \in \pi(y)> \ppp^n(x,y)} (\pi(y) - \ppp^n(x,y)) \\
& = & \sum_{y} \pi(y) \Bigl(1 - \frac{\ppp^n(x,y)}{\pi(y)} \Bigr) \\
& \le & \sep(x,n) \, .
\end{eqnarray*}

Thus  separation bounds total variation. Now we observe that it also
controls entropy decay up to a factor of $\log(1/\pi_*)$. Recall
that the relative entropy, denoted by $\D(\mu\| \nu)$, of a
distribution $\mu$ with respect to $\nu$ is defined as
\[\D(\mu \| \nu) = \sum_{x} \mu(x) \log(\mu(x)/\nu(x)),\]
where as usual $0 \log 0 = 0$, and $\mu$ is assumed to be absolutely
continuous with respect to $\nu$ (meaning,  $\nu(x)=0$ implies that
$\mu(x)=0$.) It is  well-known (see e.g. \cite{CT}) that $D(\mu\|
\nu) \le \log(1/\nu_*)$, where $\nu_* = \min_x\nu(x)$, and that
$D(\cdot \| \nu)$ is convex in the sense that, for $0 \le \alpha \le
1$, and for $\mu_1$, $\mu_2$ probability distributions (absolutely
continuous with respect to $\nu$),
\[\D(\alpha \mu_1 + (1- \alpha) \mu_2 \| \nu) \le \alpha \D(\mu_1 \| \nu) +
 (1-\alpha) \D(\mu_2\| \nu)\, .\]

\begin{proposition}
\[\D(\ppp^n(x,y) \| \pi) \le \sep(x,n) \log(1/\pi_*)\, . \]
\end{proposition}
\begin{proof}
Let $\sep(x,n) = \epsilon > 0$. Then $\ppp^n(x,y) \ge (1-\epsilon)
\pi(y)$, for all $x,y \in \xxx$. Let
\[\mu(y) := (1/\epsilon) [\ppp^n(x,y) - (1-\epsilon) \pi(y)]\, , \ \mbox{
 for } y\in \xxx\, .\]
Then $\mu$ is a probability distribution on $\xxx$ and $\ppp^n(x,\cdot) =
 (1-\epsilon) \pi + \epsilon \mu\, .$ (Note here that $\mu$ implicitly depends
on $x$.)
By the convexity mentioned above,
\[ \D(\ppp^n(x,y) \| \pi) \le \epsilon\, \D(\mu \| \pi) \le \epsilon\, \log
(1/\pi_*),\]
hence the proposition.
\end{proof}

\begin{definition}
For $0< \epsilon < 1/2$, let the entropy decay (mixing) time be
\[\tau_{{\rm ent}}(\epsilon) = \min\{n' : \ n\ge n' \ \Rightarrow \
\ \ \max_{x\in  \xxx} \ \D(\ppp^n(x,\cdot) \| \pi) \le \epsilon
\}\,.\] and similarly define the other mixing times $\tau_s, \
\tau_{\rm tv}$, and $\tau_2$ with respect to $\sep(x,n), \ d_{\rm
tv}(\ppp^n(x,\cdot) - \pi)$, and  $\|\ppp^n(x,
\cdot)/\pi(\cdot)-1\|_2$, respectively.
\end{definition}

\noindent
It then follows immediately from the above proposition,
that:
\begin{corollary}
\[\tau_{{\rm ent}}(\epsilon) \leq \tau_s(1/e) [\log \log(1/\pi_*) + \log(1/
\epsilon)]\,. \]
\end{corollary}

It is known that $\tau_s = O\bigl(\tau_{\rm tv}(\ppp) + \tau_{\rm
tv}(\ppp^*)\bigr)$, where $\ppp^*$ denotes the time-reversal of
$\ppp$.  Hence the assertion claimed at the beginning of this
section follows. Note that the lower bound below does not need
reversibility, and uses the general inequality (known as
Pinsker's) for two distributions, $\mu$ and $\nu$, one has:
\[2d^2_{\rm tv}(\mu-\nu) \le D(\mu\|\nu)\,.\]
\begin{corollary}
\label{cor:ent_tv}
 If $\ppp$ is reversible, then
\[\tau_{\rm tv}(\epsilon/2) \le \tau_{\rm ent}(\epsilon) \le C \, \tau_{\rm tv}(1/2e)\, [\log \log(1/\pi_
*) + \log(1/\epsilon)]\, ,\]
for $C>0$ an absolute constant.
\end{corollary}

\begin{remark}
Note that the above result shows that the entropy decay time is in
general closer to $\tau_{\rm tv}$ than to $\tau_2$, since there can
be a factor of $\log(1/\pi_*)$ between $\tau_{\rm tv}$ and $\tau_2$.
Similarly this indicates that in general $\rho_0$ of a reversible
Markov chain is closer to the spectral gap than it is to the
logarithmic Sobolev constant $\rho$.
\end{remark}

%It was asked in \cite{MT} whether the inverse {\it entropy
%constant} (denoted as $\rho_0$) can be used to lower bound the
%(discrete) variation mixing time of a Markov Chain. The corollary
%above shows that it is true, up to a $\log\log (1/\pi_*)$ factor,
%for reversible Markov chains.

\

We now show that the $\log\log 1/\pi_*$ gap in the corollary
cannot be improved. While a random walk on the complete graph
(with self-loops) can be shown to establish this, we proceed with
the following more robust example, which also separates various
other mixing times.

\begin{example}
Consider the following lamplighter chain on a discrete circle of
size $n$. Unlike the usual lamplighter walk, in this chain each
vertex of the circle has an $m$-state lamp for some parameter $m$.
However, every time a vertex is visited, the lamplighter
completely randomizes the lamp. The (discrete) mixing time of the
chain is still related to the time it takes to visit all vertices
of the base graph. In particular, it is easy to see that the total
variation mixing time is the expected cover time, i.e.
$\Theta(n^2)$ and is independent of $m$. From our result above, it
then follows that the entropy mixing time of this chain is
$O(n^2\log \log N)$ where $N=m^{n}$. Thus we have
\begin{equation*}
\tau_{\rm ent} = O(n^2\cdot(\log n + \log\log m)).
\end{equation*}

Suppose we start the chain at vertex $x$ of the circle and all lamps
are in state $0$. Let $A$ denote the semi-circle consisting of
vertices which are at distance $n/4$ or larger from $x$.

There exists an absolute constant $c > 0$ such that for all $a >
0$, the probability that a random walk on the circle has not
touched $A$ after $an^2$ steps is $\geq \exp(-ca)$. Thus after
$an^2$ steps, the entropy is at least $\exp(-ca)|A|\log m$ since
the entropy of the product distribution is the sum of the
component entropies and the entropy of each non-random lamp is
$\log m$. Thus in order for this chain to mix in entropy we must
have $a = \Omega(\log(n\log m))$. Hence $\tau_{\rm ent} =
\Omega(n^2\log(n\log m))$ matching the upper bound given by
Corollary~\ref{cor:ent_tv}.

Now let us look at the $L^2$ mixing time. In this case, we need to
bound the $L^2$ distance of a product space in terms of the
independent component $L^2$ distances. Since

\[\|{\mu_t}/{\pi}-1\|^2_{2,\pi} = {\rm E}_{\pi}(\mu_t/\pi)^2 -
1, \]
%the square of the $L^2$-distance is one less than the expected
%value of the square of the density function,
it follows that for $t=an^2$, we have
\begin{equation*}
1+{\rm Var}_{\pi}(\mu_t/\pi) \geq \exp(-ca) \Bigl(\prod_{i\in A}
m\Bigr) = \exp(-ca)m^{n/2}\,.
\end{equation*}
This gives $\tau_2 = \Omega(n^2 \log N)$.
\end{example}

\noindent We observe the following from the above example:
\begin{itemize}
\item If the number of states is $N=m^n$, then we have a $\Theta(\log N)$
 gap
between the variation and $L^2$ mixing times.
\item
\cite{LampLighter} shows that the relaxation time of Lamplighter
chains (with 2-state lamps) equals the maximal hitting time of the
base chain.  The $\Omega(\log N)$ gap between $\tau_{tv}$ and
$\tau_2$ also shows that in this $m$-state lamp case, the
relaxation time of the chain is $\Theta(n^2)$.
\item $\tau_{\rm ent}=\Theta(\tau_{tv}\log\log N)$ here as well.
%%This shows that the inverse entropy constant is $\tau_{tv}$.
\item Finally that we have a chain where the variation, entropy and
$L^2$-mixing times are all different orders of magnitude.
\end{itemize}

Note that in the above case the variation mixing time and the
relaxation time are of the same order of magnitude.  To separate
these, we need to separate the maximal hitting time and the
expected cover time of the underlying chain. So it is natural to
consider, the $m$-state lamplighter chain on the {\em
two-dimensional} torus of size $n$. Recall that for this case of
the torus, the maximal hitting time is $\Theta(n\log n)$ and the
expected cover time is $\Theta(n\log^2 n)$. To show that the
entropy time is still separated from the other times, observe that
there exists an $0<\alpha<1$ such that, for $a=a(n)$,
\[\Pr[\mbox{number of vertices not visited in $an \log n$ steps} \ge n^\alpha]
\ge e^{-ca},\] where $c>0$ is a constant. Now the argument follows as
in the one-dimensional case above, and gives $\tau_{\rm ent} =
\Omega(n\log(n \log m))$.

\section{Questions}

Suppose $\ppp,\qqq$ are Markov Chains on state space $\xxx,\yyy$
respectively. The Lamplighter chain $\qqq \wr \ppp$ has state space
$\yyy^\xxx \times \xxx$, i.e. a configuration of lamps $f$ together
with the position of the lamplighter. \cite{LampLighter}
considered $\yyy=\zzz_2$ and $\qqq$ on $\yyy$ which completely
randomizes the lamp in one step.

In \cite{ProdChains} it is shown that the $L_2$-mixing time of the
lamplighter chain on $\qqq \wr \ppp$ is related to the following
generalization of the moment generating function considered in this
paper.

\begin{definition}
Let $\ppp$ be a Markov Chain on state space $\xxx$ and let $\gamma
> 0$. For $S > 0$, let $\rv{Z}_S^{\gamma}$ denote the number of
states that have been visited fewer than $\gamma$ times up till
time $S$ and let $\zeta_S^{\gamma}(\theta) =
\EXP[\theta^{\rv{Z}_S^{\gamma}}].$
\end{definition}

Quantities similar to $\rv{Z}_S^{\gamma}$, with first moment
computations, were considered in \cite{Blanket} and
\cite{KahnBlanket}. The first paper to consider moment generating
functions of Markov chain related quantities seems to be
\cite{LampLighter}.

\begin{question}
For $\theta > 0,\delta > 0$, find bounds on
\begin{equation*}
F(\ppp,\theta,\gamma,\delta) = \inf_S \{ S: \zeta_S^\gamma(\theta)
\leq 1+\delta \}\,.
\end{equation*}
For estimating the mixing time of $\qqq \wr \ppp$, where $\qqq$ is
a Markov chain on $\yyy$, the quantity of interest is
$F(\ppp,|\yyy|,\Tau_2(\qqq,\epsilon/|\xxx|),\epsilon)$.
\end{question}

\begin{itemize}
\item If $\gamma \geq |\xxx|\log\theta$, then it is enough to take
$S=O(\gamma|\xxx|)$.
\item If $\ppp$ mixes in one step then it reduces to the coupon
collector problem and hence for any $\gamma \geq 0$, it is enough to
take $S=O((\gamma+\log|\xxx|)|\xxx|)$.
\item In general,
$O((\gamma+\log|\xxx|)|\xxx|T_{\rm tv})$ is enough, where $T_{\rm
tv}$ is the variation mixing time of $\ppp$.
\end{itemize}

\

\noindent \textbf{Conjecture 1}. Let $\ppp$ be reversible with
uniform stationary distribution and maximal hitting time
$H=O(|\xxx|)$. Then
\begin{equation}
F(\ppp,\theta,\gamma,\delta) \leq C\cdot |\xxx|
(\gamma+T_{rel}+\log|\xxx|)\, ,
\end{equation}
for some absolute constant $C$.

\

\noindent {\bf Acknowledgment}. We thank Yuval Peres for general
encouragement and for help with examples in Section 3.

%%\bibliography{touching}

\begin{thebibliography}{10}
\providecommand{\natexlab}[1]{#1}
\providecommand{\url}[1]{\texttt{#1}} \expandafter\ifx\csname
urlstyle\endcsname\relax
  \providecommand{\doi}[1]{doi: #1}\else
  \providecommand{\doi}{doi: \begingroup \urlstyle{rm}\Url}\fi

\bibitem[Aldous and Diaconis(1987)]{AD}
D.~Aldous and P.~Diaconis.
\newblock Strong uniform times and finite random walks.
\newblock \emph{Advances in Appl. Math.}, 8:\penalty0 69--97, 1987.

\bibitem[Aldous and Fill()]{AF}
D.~Aldous and J.A. Fill.
\newblock Preliminary version of a book on finite markov chains.
\newblock URL \url{http://www.stat.berkeley.edu/~aldous/RWG/book.html}.
\newblock In preparation.

\bibitem[Aldous(1982)]{Aldous}
David Aldous.
\newblock Some inequalities for reversible {M}arkov chains.
\newblock \emph{J. London Math. Soc. (2)}, 25\penalty0 (3):\penalty0 564--576,
  1982.

\bibitem[Cover and Thomas(1991)]{CT}
T.M. Cover and J.A. Thomas.
\newblock \emph{Elements of Information Theory}.
\newblock John-Wiley \& Sons, Inc., 1991.

\bibitem[Ganapathy(2006)]{ProdChains}
Murali~K. Ganapathy.
\newblock \emph{Robust Mixing}.
\newblock PhD thesis, University of Chicago, August 2006.

\bibitem[Horn and Johnson(1991)]{HJ}
R.A. Horn and C.R. Johnson.
\newblock \emph{Topics in Matrix Analysis}.
\newblock Cambridge University Press, 1991.

\bibitem[Kahn et~al.(2000)Kahn, Kim, Lov\'{a}sz, and Vu]{KahnBlanket}
J.~Kahn, J.H. Kim, L.~Lov\'{a}sz, and V.H. Vu.
\newblock The cover time, the blanket time, and the matthews bound.
\newblock In \emph{41st Annual Symposium on Foundations of Computer Science},
  pages 467--475, 2000.

\bibitem[Montenegro and Tetali(2006)]{MT}
R.~Montenegro and P.~Tetali.
\newblock Mathematical aspects of mixing times of markov chains.
\newblock \emph{Current and Future Trends in Theoretical Computer Science}, May
  2006.

\bibitem[Peres and Revelle(2004)]{LampLighter}
Y.~Peres and D.~Revelle.
\newblock Mixing times for random walks on finite lamplighter groups.
\newblock \emph{Electronic Journal of Probability}, 9:\penalty0 825--845, 2004.
\newblock URL
  \url{http://www.math.washington.edu/~ejpecp/EjpVol9/paper26.abs.html}.

\bibitem[Winkler and Zuckerman(1996)]{Blanket}
P.~Winkler and D.~Zuckerman.
\newblock Multiple cover time.
\newblock \emph{Random Structures and Algorithms}, 9\penalty0 (4):\penalty0
  403--411, 1996.
\newblock URL \url{citeseer.ist.psu.edu/139.html}.

\end{thebibliography}

\end{document}